\documentclass[12pt]{amsart}
\usepackage{amsfonts, amsbsy, amsmath, amssymb}
\usepackage{multicol} 
\usepackage{algorithm}
\usepackage{algpseudocode}
\usepackage{tikz}
\usepackage{array}
\usepackage{longtable}
\usepackage{xltabular}
\usepackage{tabularx}
\usepackage{rotating}
\usepackage{mathtools}
\usepackage{url}
\usepackage{lscape}

\usetikzlibrary{knots}

\usepackage[lite]{amsrefs}

\newtheorem{thm}{Theorem}[section]

\newtheorem{exmp}[thm]{Example}

\newtheorem{conj}[thm]{Conjecture}

\newtheorem{thm-con}[thm]{Theorem-Conjecture}
\numberwithin{equation}{section}

\theoremstyle{definition}
\newtheorem{definition}[thm]{Definition}
\newtheorem{defn}[thm]{Definition}

\title[Isomorphism Classes of Connected Shelves]
{Classification of Connected Shelves}

\begin{document}

	\author[Elhamdadi]{Mohamed Elhamdadi} 
	\address{Department of Mathematics and Statistics, University of South Florida, Tampa, FL, USA}
	\email{emohamed@usf.edu}

 	\author[Fernando]{Neranga Fernando} 
	\address{Department of Mathematics and Computer Science,
College of the Holy Cross, Worcester, MA, USA}
	\email{nfernand@holycross.edu}

 	\author[Goonewardena]{Mathew Goonewardena} 
	\address{Ericsson, Montreal, QC, Canada}
	\email{mathew.goonewardena@ericsson.com}

\maketitle 

\begin{abstract}
We investigate finite right-distributive binary algebraic structures called shelves. We first use symbolic computations with Python to classify (up to isomorphism) all connected shelves with order less than six. We explore the group structure generated by the rows of \textit{latin} shelves. We also define two-variable shelf polynomial by analogy with the quandle polynomial and then state a conjecture about connected idempotent shelves. 
 \end{abstract}
 \tableofcontents

\section{Introduction}

Shelves are sets with binary operations satisfying self-distributivity $$(x*y)*z=(x*z)*(y*z).$$  These algebraic structures are derived from the axiomatization of Reidemeister move III in classical knot theory.  
By adding the condition corresponding to Reidemeister move II, one gets the notion of a rack.  Racks have been used to obtain invariants of framed knots \cite{FR}.  Framed knots can be visualized as closed loops of knotted flat ribbons. Framed knots generated lot of interests mainly because of the critical role they play in low-dimensional topology \cite{EH, FR}.  By adding the condition corresponding to Reidemeister move I to the definition of a rack, one obtains the notion of a quandle.  Quandles have been studied extensively and are used to obtain invariants of knots and links \cite{Joyce, Matveev, EN}. In \cite{CMP}, associative shelves have been investigated and it was shown that unital shelves are associative. The authors also investigated one-term and two-term homology groups of some associative self-distributive algebraic structures in \cite{CMP}.  
In \cite{CCES}, self-distributivity was investigated in a unified manner via a categorical technique called internalization, and a cohomology theory was developed  and explicit relations to rack and Lie algebra cohomology theories were given. Self-distributivity also provides solutions of the Yang-Baxter equation \cite{CES}. There has been other investigations in relations to many other areas of mathematics such as ring theory \cite{EFT, BPS1, BPS2, ENSS, ENS}, quasigroups and Moufang loops \cite{E}, representation theory \cite{EM} and singular knot theory \cite{CEHN, BEHY, CCE}.  

 A \emph{monounary} algebra is an algebra with one unary operation, usually denoted as a pair $(A,f)$, where $A$ is a nonempty set with a map $f:A \rightarrow A$.  In \cite{Jezek}, the author considered finite monounary algebras and proved that every monounary algebra has at least one left-distributive extension $E$ such that $x*y=z$ in $A$ means $x*y=z$ in $E$. A left-distributive extension of $A$ is simply a left-distributive groupoid. The author also gave an enumeration algorithm which computes the numbers of all left-distributive groupoids and their isomorphism types on a given set of cardinality less or equal to six. It is worth noting that the left-groupoids discussed in \cite{Jezek} are exactly the shelves that are under consideration in our paper, but with right-distributive property. However, to the best of our knowledge, a classification of connected shelves up to isomorphism has not appeared in the literature. 

A Laver table is a free shelf generated by one element.  Precisely, Laver proved that for every $N\geq 1$, there exists a unique binary operation $\ast$ on the set $X=\{1, 2, \dots, N\}$ that, for all $x,y$, satisfies 
\begin{center} 
$x*1=x+1\pmod{N},$
$x\ast (y\ast 1)=(x\ast y)\ast (x\ast 1)$
\end{center} 
Then the binary operation $\ast$ is left-distributive if and only if $N$ is a power of 2. Laver tables were introduced in 1995 by Richard Laver while investigating self-embedding in set theory.  Recently they have been investigated from the topology point of view, see for example \cite{Deh}, which discusses the use of Laver tables in low-dimensional topology, and \cite{Dehornoy-Lebed}, which classifies $2$- and $3$-cocycles on Laver's tables.   

In this article, we use symbolic computations with Python to obtain the list of shelves of order less than six. We also define two-variable shelf polynomial by analogy with the quandle polynomial and then state a conjecture about connected idempotent shelves. The organization of this paper is as follows: Section~\ref{review} reviews the basics of shelves in general and connected shelves in particular.  In Section~\ref{algorithm}, we describe the algorithm which allows the classification of connected shelves of order less than six. In Section~\ref{Connected}, we study connected shelves and give the number of connected shelves of order up to 6. We also present a conjecture regarding connected shelves. A {\it latin} shelf is a shelf whose rows are permutations. We also study latin shelves in Section~\ref{Connected}, especially the group structure generated by rows of latin shelves. Section~\ref{polynomial} introduces shelf polynomial, and we end the section with a conjecture strongly supported by our computational results. In Appendix A, we give the list of all connected shelves of order less than six. 

\section{Review of Shelves}\label{review}
In this section, we review the basics of shelves, give some examples and introduce the notion of connected shelves.

\begin{definition}\label{Shelf}
  A \emph{shelf} $(S,*)$ is a set $S$ with a binary operation $*$ that is right-distributive.  Precisely, for all $x,y,z \in S$, we have:
\[
(x*y)*z=(x*z)*(y*z).
\]  
\end{definition}
A shelf homomorphism between two shelves $(S,*_1)$ and $(S',*_2)$ is a map $f:S \rightarrow S'$ such that for all $x,y \in S$ we have $f(x*_1y)=f(x)*_2f(y)$.  A shelf \emph{isomorphism} is a bijective shelf homomorphism.   Two shelves are isomorphic if there is a
shelf isomorphism between them.

Typical examples of shelves include:
\begin{itemize}\item
The set $\mathbb{Z}_n$ of integers modulo $n$ with $x*y=\alpha x +\beta y$, such that $\beta(\alpha+\beta -1) =0 $ in $\mathbb{Z}_n$.  For example $\mathbb{Z}_{10}$ with $x*y=2x+5y$.
    \item 
Any group $G$ with conjugation $x*y=yxy^{-1}$. 
\item
Laver tables: $\mathbb{Z}_{2^n}$ with $1*x=x+1$ and $(1*x)*y=(1*y)*(x*y)$.  (Notice that we are using right-distributivity instead of left-distributivity).
\end{itemize}
Let $(S,*)$ be a shelf.  The right multiplication $R_x$ by an element $x$ of $S$ is the map $R_x: S \rightarrow S$ such that $R_x(y)=y*x$.  The monoid generated by $R_x$ is denoted by $M_S$.  This monoid
acts naturally on $S$.  When this action is transitive we say that the shelf $S$ is \emph{connected}.  In other words, if for all $x,y\in S$, there exists $\phi \in M_S$ such that $\phi(x)=y$, then we say that the shelf $S$ is connected. 

A {\it unital} shelf is a shelf $S$ with an identity element. That is, there exists an element 1 in $S$ such that $1\ast x = x = x\ast 1$ for all $x\in S$. A {\it spindle} is a shelf whose elements are idempotent. In \cite{CMP}, the authors developed a theory of associative shelves, and studied their one-term and two-term homology groups. They also presented the number of unital shelves of order up to 4 in Table 3. The following table is an update of it. 

\vskip 0.3in

\begin{center} 
 \begin{tabular}{|c|c|} 
\hline
 \hfil $n$ & number of unital shelves\\
\hline
\hfil 1\hfil  &\hfil 1 \\
    \hfil 2 & \hfil  1\\
    \hfil 3 & \hfil  2\\ 
    \hfil 4 & \hfil  6\\
    \hfil 5 & \hfil  23 \\
   \hline
\end{tabular}
\end{center} 

\section{An Algorithm to Generate Shelves} \label{algorithm}

A matrix $M \in M_n(\mathbb{Z}_n)$ is a shelf matrix if it satisfies right distributivity [ref to definition]. A brute-force algorithm that enumerates all matrices in $M_n(\mathbb{Z}_n)$ and check each candidate for right distirbutivity does not scale for orders higher than four, e.g., order five requires checking $5^{25}$ candidate matrices. This section presents a computationally efficient algorithm to generate shelves without enumerating all of $M_n(\mathbb{Z}_n)$. The proposed algorithm is also parallelizable.

The right distributivity property requires to check $n^3$ conditions (all possible orderings of three elements drawn with replacement from the set). Let us denote this set of conditions by $\mathcal{C}_n$ and one of those conditions by a tuple of length three $(x,y,z)$, where $x,y,z \in \mathbb{Z}_n.$ For example, consider order three where $\mathbb{Z}_3 = \{0,1,2\}$ and $\mathcal{C}_n \coloneqq \{(x,y,z) \mid x,y,z \in \mathbb{Z}_n\}$ is the set of conditions and $(0,1,2)$ denotes the particular right distirbutivity condition $(0\ast1)\ast 2 = (0\ast 2) \ast (1\ast 2).$

The algorithm starts with an empty $n \times n$ matrix. Then it takes one of the conditions from $\mathcal{C}_n$ and  generates all the possible partially filled matrices that satisfy the condition. We call these candidate matrices. For example consider $(0,0,0)$ condition in order three. Let us represent the empty matrix as $$[[ -1, -1, -1], [-1, -1, -1], [-1, -1, -1]],$$ where the rows are unstacked into a single row to preserve space and $-1$ is used to represents the empty value in the computer program (any value outside of $\mathbb{Z}_3$ will do). From this empty matrix generating all partially filled matrices that satisfy condition $(0,0,0)$ gives the following seven candidates:

\begin{enumerate} 
\item [ (i)] $[[0, -1, -1], [-1, -1, -1], [-1, -1, -1]]$
\item [(ii)] $[[1, -1, -1], [\,\,\,\,0, \,\,\,\,\,0, -1], [-1, -1, -1]]$
\item [(iii)]$[[1, -1, -1],[\,\,\,\,1, \,\,\,\,\,1, -1], [-1, -1, -1]]$
\item [(iv)] $[[1, -1, -1], [\,\,\,\,2, \,\,\,\,\,2, -1], [-1, -1, -1]]$
\item [ (v)] $[[2, -1, -1], [-1, -1, -1], [\,\,\,\,\,0, -1, \,\,\,\,0]]$
\item [(vi)] $[[2, -1, -1], [-1, -1, -1], [\,\,\,\,\,1, -1, \,\,\,\,1]]$
\item [(vii)] $[[2, -1, -1], [-1, -1, -1], [\,\,\,\,\,2, -1, \,\,\,\,2]]$ 
\end{enumerate} 

The generation technique of the candidates is as follows. In Python or C programming array indices start at $0$. Therefore, the operation $(x\ast y)$, where $x,y \in \mathbb{Z}_n,$ is identical to reading the value at $(x,y)$ indices of the 2D array containing the matrix. Assume we are given a matrix $M$ (partially or fully empty) and the condition $(x,y,z)$. First check if the location $(x,y)$ is empty in $M$ and if so it can be filled with $n$ possible values. This is done in a loop. Enter the loop and fill with one value. Now check if $(x,z)$ is empty and if so it too can be filled with $n$ possible values in a similar loop. Do the same for $(y,z)$. Then in the final stage check if either one of $((x,y),z)$ or $((x,z),(y,z))$ is empty. If both are empty (or they point to same location in the matrix and that is empty) there are $n$ possible values to consider, each of which satisfies the condition. If just one is empty then assign the value of other to it. If both are filled with the same value then the condition is satisfied. If they are filled but with different values then the condition is violated. The generation process is computationally intensive with a series of nested loops. Therefore, an efficient algorithm should minimize the total number of false candidates that it generates while searching for shelves. This beings to the next important step of the algorithm.

For each of these candidates check if at least one of the remaining conditions is violated. For example consider the condition $(0,0,1)$ and first candidate $$[[ \,\,0, -1, -1], [-1, -1, -1], [-1, -1, -1]]$$ generated above. To check both LHS and RHS of $(0,0,1)$ we need the value at $(0, 1)$. However, $(0,1) = -1$, as the matrix location $(0, 1)$ has not yet been filled by a value from $\mathbb{Z}_3$. Therefore, one cannot make a conclusion as to if this conditions is violated.  Likewise check this matrix for all of the remaining conditions for violations. This checking is important to reduce the search space in the latter steps. Those that do not violate any of the remaining conditions become candidates for the next step. Now in the next step the algorithm takes each of these candidates and the next condition and generates all possible candidates that do not violate the remaining conditions. Thus, in each step the candidates satisfy all the conditions that have been applied so far to generate them and they also do not violate the remaining conditions. This process can be visualized as a tree with the empty matrix at the root.

The above described steps are succinctly presented in Algorithm \ref{alg:shelves} in the form of a recursive depth first search. The function \texttt{getNextCondition()}, takes as input the current condition (\texttt{None} for the first call) and returns the next condition to apply. The function \texttt{generateCandidates()} takes as input a candidate matrix, and a condition and then derives all candidates that satisfy the given condition. The function \texttt{notViolateRestOfTheConditions()} takes as input a matrix and the current condition and returns True if the matrix does not violate any of the remaining conditions. This function can be modified to verify additional conditions when looking for specific kinds of shelves, such as connected or unitial. The algorithm is started by calling the function \texttt{getShelves()} with the empty initial matrix $M_\text{init}$ and the fist condition. Then \texttt{getShelves()} proceeds recursively in a depth first manner collecting the matrices that satisfy all the conditions of right distributivity, which are shelves.

A Python implementation of Algorithm \ref{alg:shelves}, a list of shelves of order less than 6 (up to isomorphism), and a list of connected shelves of order less than seven (up to isomorphism) are available in \cite{Aloglink}. This Python implementation also parallilizes the algorithm for multiprocessing. The parallization is achieved by first applying two of the $n^3$ conditions in a breadth first search, thus generating all the candidates that satisfy these two conditions and that do not violate remaining conditions and then parallelly applying the depth first search on these candidates. The order in which the conditions are applied does not affect the final result, however, through experimentation we found that the order very much affects the time to completion. Optimal ordering is left for future research.


\begin{algorithm}
\algrenewcommand\textproc{}
\caption{Algorithm to construct shelves}\label{alg:cap}
\label{alg:shelves}
\begin{algorithmic}
\State $n \gets \text{order}$ \Comment{$n \geq 2$ }
\State $M_{\text{init}} \gets \text{empty matrix of size } n \times n$
\State $c_1 \gets $ \Call{getNextCondition}{None}
\State $L \gets [ ]$ \Comment{list to store shelves}
\Procedure{getShelves}{$M$,$c_{\text{in}}$} \Comment{recursively applies conditions}
  \State $c \gets$ \Call{getNextCondition}{$c_{\text{in}}$}
  \If{\text{all conditions applied} $\And${$M$\text{ not in }$L$}}
  
    \State append $M$ to $L$
    
  \Else
    \For{$m$ in \Call{generateCandidates}{$M$, $c_{\text{in}}$}}
        \If{\Call{notViolateRestOfTheConditions}{$m$, $c_{\text{in}}$}}
        \State \Call{getShelves}{$m$, $c$}
        \EndIf
    \EndFor
  \EndIf
\EndProcedure
\State \Call{getShelves}{$M_{\text{init}}$, $c_1$} \Comment{ call the procedure}
\State \Comment{when the procedure exits $L$ contains all the shelves of order $n$}
\end{algorithmic}
\end{algorithm}


\newpage

\section{Connected Shelves} \label{Connected}

In this section, we classify connected shelves of order up to 5. 

\begin{defn}
A connected shelf $(X,*)$ is a shelf such that for all $x,y\in X$, there exists a finite number of elements $x_1,x_2,\ldots , x_m$ such that $$y=((((x\ast x_1)\ast x_2)\ast x_3)\ldots )\ast x_m$$
\end{defn} 

The definition of a connected shelf simply means that we can go from one element to any other element by finite number of steps. In other words, the orbit of each element must equal the shelf itself. 

\vskip 0.2in 

\begin{tabular}{|c|c|c|c|} 
\hline
 \hfil $n$ & \# of connected shelves & \# of connected racks & \# of connected quandles \\
\hline
\hfil 1\hfil  &\hfil 1 &\hfil 1 & 1  \\
    \hfil 2 & \hfil 2 &\hfil 1 & 0\\
    \hfil 3 & \hfil 5 &\hfil 2 & 1 \\
    \hfil 4 & \hfil 18 &\hfil 2 & 1 \\
    \hfil 5 & \hfil  165 &\hfil 4 & 3 \\
    \hfil 6 & \hfil 3987 &\hfil 4 & 2 \\
   \hline
\end{tabular}

\vskip 0.2in 

We refer the reader to Appendix A for the complete list of connected shelves of order less than or equal to 5, up to isomorphism. 

\vskip 0.2in

Our computer search results support the following conjecture for order up to 5. 

\begin{conj}\label{C1}
Let S be a shelf.  Then there exists a shortest cycle that covers all the elements in $S$ exactly once if and only if S is connected.

\end{conj} 

A cycle from 0 to 0 that covers all elements in the shelf may contain an element more than once. For example, consider the following connected shelf of order 4. 

\begin{center}
\begin{tabular}{ r| c c c c}
* & 0 & 1 & 2 & 3\\
\hline
  0 & 0 & 1 & 1 & 3\\
  1 & 0 & 1 & 2 & 0 \\
  2 & 0 & 1 & 2 & 0 \\
  3 & 0 & 1 & 1 & 3 \\
\end{tabular} 
\end{center}

A cycle from 0 to 0 in this shelf is 

$$0\rightarrow 1\rightarrow 2\rightarrow 1\rightarrow 3\rightarrow 0.$$

However the shortest cycle from 0 to 0 is 

$$0\rightarrow 3\rightarrow 1\rightarrow 2\rightarrow 0.$$

\subsection{Latin Shelves} 

In this subsection, we explore the groups generated by rows in Latin shelves. Recall the left multiplication map. For each $x\in S$, the left multiplication by $x$ is the map denoted by $$L_x\,:\,S\rightarrow S$$ and given by $$L_x(y)\,:=\,x\ast y.$$

\begin{defn}
A shelf is {\it Latin} or {\it strongly connected} if the shelf operation is left-invertible. This means the rows of a Latin shelf are permutations of $\{1,2,\ldots , n\}$. 
\end{defn} 

Clearly, Latin shelves are connected, and a Latin quandle is always a Latin shelf. 

Let $z\in S$, where $S$ is a Latin shelf. 

$$L_{x\ast y}(z)=(x\ast y)\ast z=(x\ast z)\ast (y\ast z)=L_x(z)\ast L_y(z)$$

In a Latin shelf, each $L_x$, where $x\in S$, is a bijection. 

Let $X$ be the set of all $L_x$, i.e. 

$$X=\{L_x\,:\,x\in S\}.$$

Clearly, we have $X\subset S_n$. 

Let $t\in S$. Define an operation $\rhd$ on the set $X$ as follows:

$$(L_x\rhd L_y)(t)\,:=\,L_x(t)\ast L_y(t).$$

Then $$(L_x\rhd L_y)(t)=L_{x\ast y}(t).$$

Now consider 

\[
\begin{split}
((L_x\rhd L_y)\rhd L_z)(t)&=(L_x\rhd L_y)(t)\ast L_z(t)\cr
&=(L_x(t)\ast L_y(t))\ast L_z(t)\cr
&=(L_x(t)\ast L_z(t))\ast (L_y(t)\ast L_z(t))\cr
&=(L_x\rhd L_z)(t) \ast (L_y\rhd L_z)(t)\cr
&=L_{x\ast z}(t)\ast L_{y\ast z}(t)\cr
&=(L_{x\ast z}\rhd L_{y\ast z})(t)\cr
&=((L_x\rhd L_z)\rhd (L_y\rhd L_z))(t)
\end{split}
\]

Since $t$ is arbitrary, we have

$$(L_x\rhd L_y)\rhd L_z=(L_x\rhd L_z)\rhd (L_y\rhd L_z),$$

i.e.

$$L_{x\ast y}\rhd L_z=L_{x\ast z}\rhd L_{y\ast z}.$$

Thus, $(X,\rhd)$ is a shelf. 

\vskip 0.1in

Note that $L_x\rhd L_x=L_{x\ast x}$. So, if the elements in the shelf $S$ are idempotent, so are the elements in $X$. 

Now we show that the elements in $X$ do not satisfy the second axiom of a quandle. 

Let $t\in S$. 

\[
\begin{split}
((L_x\rhd L_y)\rhd L_y)(t)&=(L_x\rhd L_y)(t)\ast L_y(t)\cr
&=(L_x(t)\ast L_y(t))\ast L_y(t)\cr
&=(L_x(t)\ast L_y(t))\ast (L_y(t)\ast L_y(t))\cr
&=(L_x\rhd L_y)(t) \ast (L_y\rhd L_y)(t)\cr
&=L_{x\ast y}(t)\ast L_{y\ast y}(t)\cr
&=L_{(x\ast y)\ast (y\ast y)}(t)\neq L_x(t)
\end{split}
\]

For each $x\in S$, define a map
$$\phi\,:\,S\rightarrow X$$

given by

$$x\mapsto L_x$$

The map $\phi$ is a homomorphism because

$$\phi(x\ast y)=L_{x\ast y}=L_x\rhd L_y=\phi(x)\rhd \phi(y)$$

In fact, it is an epimorphism. 

Now let $G$ be the group generated by $L_x$, where $x\in S$, i.e. $G=\langle L_x\,:\,x\in S\rangle$. Since we are only considering finite shelves, $G$ is a finitely generated group. 

Define an operation on $G$ as 

$$L_x\diamond L_y\,:=\,L_y^{-1}L_xL_y$$

Then we have 

\[
\begin{split}
(L_x\diamond L_y)\diamond L_z&=L_z^{-1}(L_x\diamond L_y)L_z\cr
&=L_z^{-1}(L_y^{-1}L_xL_y)L_z\cr
&=(L_yL_z)^{-1}L_x(L_yL_z)
\end{split}
\]

and

\[
\begin{split}
(L_x\diamond L_z)\diamond (L_y\diamond L_z)&=(L_z^{-1}L_xL_z)\diamond (L_z^{-1}L_yL_z)\cr
&=(L_z^{-1}L_yL_z)^{-1}(L_z^{-1}L_xL_z)(L_z^{-1}L_yL_z)\cr
&=(L_z^{-1}L_y^{-1}L_z)(L_z^{-1}L_xL_z)(L_z^{-1}L_yL_z)\cr
&=(L_z^{-1}L_y^{-1})L_x(L_yL_z)\cr
&=(L_yL_z)^{-1}L_x(L_yL_z),
\end{split}
\]

which imply

$$(L_x\diamond L_y)\diamond L_z=(L_x\diamond L_z)\diamond (L_y\diamond L_z)$$

Hence the conjugation in $G$ turns it into a shelf. 

We are now interested in seeing whether the conjugation in $G$ satisfies Axiom 1 and Axiom 2 of a quandle. Consider $L_x\diamond L_x$. 

$$L_x\diamond L_x=L_x^{-1}L_xL_x=L_x$$

Thus the elements in $G$ are idempotent. 

Let's consider Axiom 2. We show that under certain conditions on left multiplication, elements in $G$ satisfy Axiom 2. In other words, under certain conditions on left multiplication, $G$ is a quandle. 

\[
\begin{split}
(L_x\diamond L_y)\diamond L_y&=L_y^{-1}(L_x\diamond L_y)L_y\cr
&=L_y^{-1}(L_y^{-1}L_xL_y)L_y\cr
&=(L_y^{-1})^{2}L_x(L_y)^2\cr
&=(L_y^{2})^{-1}L_x(L_y^2)
\end{split}
\]

If rows are either 2-cycles or the identity permutation, then $(L_x\diamond L_y)\diamond L_y=L_x$, and thus $(G,\diamond)$ is a quandle. If cycles are disjoint, $(L_x\diamond L_y)\diamond L_y=L_x$, and thus $(G,\diamond)$ is still a quandle. 

Since the rows of Latin shelves are permutations, we are curious about the groups generated by them. We are interested in knowing whether the group structure can be used to distinguish shelves in different isomorphism classes. In the following table, we present the groups generated by rows of Latin shelves of order up to 5. In the table, $S_{n,k}$ denotes the $k$th shelf of order $n$ listed in Appendix A. The shelf $S_{3,3}$ denotes the Latin shelf \# $3$ listed under order 3 in Appendix A. We first list disjoint cycle notation for the rows of a Latin shelf, and then we list the group generated by them. The cycle $({\rm id})$ denotes the identity permutation. 

Consider the following example: $$S_{3,3}=[[0,2,1],[2,1,0],[1,0,2]]$$
The first, second and third rows are two cycles $(12),(02),(01)$, respectively. Thus the group generated by rows (left multiplications) is 
$$G=\langle (12),(02),(01)\rangle = S_3\cong D_3$$
  
In fact, this is the only connected Latin quandle of order 3. The abbreviations LQ stand for Latin Quandle. 

\vskip 0.2in

\begin{center}
\begin{scriptsize}
 \begin{tabular}{|c|c|c|c|} 
\hline
 \hfil Latin Shelf & Disjoint Cycle Notation for the Rows of the Latin Shelf & $G$ \\
\hline
 \hfil $S_{2,1}$ &\hfil $({\rm id})$,$({\rm id})$ & \{\rm id\}  \\
 \hfill                 &\hfill           &          \\
 \hfil $S_{3,1}$ &\hfil $({\rm id})$,$({\rm id})$,$({\rm id})$ & \{\rm id\}\\
 \hfil $S_{3,2}$ &\hfil $({\rm id})$,$({\rm id})$,(01) & $\mathbb{Z}_2$ \\
 \hfil $S_{3,3(LQ)} $ &\hfil (12),(02),(01) & $D_3$ \\
 \hfill                 &\hfill           &          \\
 \hfil $S_{4,1}$ &\hfil $({\rm id})$,$({\rm id})$,$({\rm id})$,$({\rm id})$ & \{\rm id\} \\
 \hfil $S_{4,2(LQ)}$ &\hfil (132),(023),(031),(012) & $A_4$\\
 \hfil $S_{4,3}$ &\hfil (23),(23),(01),(01) & $\mathbb{Z}_2 \times \mathbb{Z}_2$  \\
 \hfil $S_{4,4}$ &\hfil $({\rm id})$,$({\rm id})$,(01),(01) & $\mathbb{Z}_2$ \\
 \hfil $S_{4,7}$ &\hfil $({\rm id})$,$({\rm id})$,$({\rm id})$,(021) & $\mathbb{Z}_3$ \\
 \hfil $S_{4,13}$ &\hfil $({\rm id})$,$({\rm id})$,$({\rm id})$,(12) & $\mathbb{Z}_2$ \\
  \hfill                 &\hfill           &          \\
  \hfil $S_{5,126}$ &\hfil $({\rm id})$,$({\rm id})$,$({\rm id})$,$({\rm id})$,$({\rm id})$ & \{\rm id\} \\
  \hfil $S_{5,127}$ &\hfil $({\rm id})$,$({\rm id})$,$({\rm id})$,$({\rm id})$,(23) & $\mathbb{Z}_2$ \\
  \hfil $S_{5,129}$ &\hfil $({\rm id})$,$({\rm id})$,$({\rm id})$,$({\rm id})$,(132) & $\mathbb{Z}_3$ \\
  \hfil $S_{5,133}$ &\hfil $({\rm id})$,$({\rm id})$,$({\rm id})$,$({\rm id})$,(01)(23) & $\mathbb{Z}_2$ \\
  \hfil $S_{5,135}$ &\hfil $({\rm id})$,$({\rm id})$,$({\rm id})$,$({\rm id})$,(0321) & $\mathbb{Z}_4$ \\
  \hfil $S_{5,136}$ &\hfil $({\rm id})$,$({\rm id})$,$({\rm id})$,(12),(12) & $\mathbb{Z}_2$ \\
  \hfil $S_{5,139}$ &\hfil $({\rm id})$,$({\rm id})$,$({\rm id})$,(12),(01) & $D_3$ \\
  \hfil $S_{5,140}$ &\hfil $({\rm id})$,$({\rm id})$,$({\rm id})$,(12),(021) & $D_3$ \\
  \hfil $S_{5,150}$ &\hfil $({\rm id})$,$({\rm id})$,$({\rm id})$,(021),(021) & $\mathbb{Z}_3$ \\
  \hfil $S_{5,151}$ &\hfil $({\rm id})$,$({\rm id})$,$({\rm id})$,(021),(012) & $\mathbb{Z}_3$ \\
  \hfil $S_{5,152}$ &\hfil $({\rm id})$,$({\rm id})$,(34),(01),(01) & $\mathbb{Z}_2 \times \mathbb{Z}_2$  \\
  \hfil $S_{5,153}$ &\hfil $({\rm id})$,$({\rm id})$,(01),(01),(01) & $\mathbb{Z}_2$ \\
  \hfil $S_{5,154}$ &\hfil $({\rm id})$,$({\rm id})$,(01),(01),(01)(23) & $\mathbb{Z}_2 \times \mathbb{Z}_2$ \\
  \hfil $S_{5,155}$ &\hfil $({\rm id})$,(34),(34),(12),(12) & $\mathbb{Z}_2 \times \mathbb{Z}_2$  \\
  \hfil $S_{5,156}$ &\hfil (34),(34),(34),(12),(12) & $\mathbb{Z}_2 \times \mathbb{Z}_2$ \\
  \hfil $S_{5,158}$ &\hfil (34),(34),(34),(021),(021) &  $\mathbb{Z}_2 \times \mathbb{Z}_3$\\
  \hfil $S_{5,159}$ &\hfil (34),(34),(01)(34),(01),(01) & $\mathbb{Z}_2 \times \mathbb{Z}_2$\\
  \hfil $S_{5,162(LQ)}$ &\hfil (12)(34),(03)(24),(13)(04),(02)(14),(01)(23) & $D_5$\\
  \hfil $S_{5,163(LQ)}$ &\hfil (1432),(0342),(0413),(0124),(0231) & $GA(1,5)$\\
  \hfil $S_{5,164(LQ)}$ &\hfil (1432),(0423),(0134),(0241),(0312) & $GA(1,5)$\\
   \hline
\end{tabular}
\end{scriptsize}
\end{center}

\vskip 0.2in

We are able to distinguish many Latin shelves that belong to different isomorphism classes, but there are quite a few cases starting from order 4 in which shelves that belong to different isomorphism classes have the same group structure. For example, in order 4, the rows of Latin shelves $S_{4,4}$ and $S_{4,13}$ have the group structure $\mathbb{Z}_2$, and in order 5, the rows of Latin shelves $S_{5,152}$, $S_{5,154}$, $S_{5,155}$, $S_{5,156}$ and $S_{5,159}$ generate the group $\mathbb{Z}_2\times \mathbb{Z}_2$. 

It also seems like the number of such cases increases as the order gets higher. However, we believe that the group structure of the rows of Latin shelves would be an interesting topic for further research, especially as the order gets higher. 

\section{Shelf polynomial} \label{polynomial}

In this section, we study shelf polynomials. 

\begin{defn} 
Let $S$ be a finite shelf. For each element $x\in S$, let $r(x)$ be the number of elements of $S$ which act trivially on $x$, i.e. the set 
$$r(x)=|\{y\in S\,|\,x\ast y= x\}|$$

and let $c(x)$ be the number of elements of $S$ on which $x$ acts trivially, i.e. the set 
$$c(x)=|\{y\in S\,|\,y\ast x= y\}|$$

In terms of the shelf's operation table, $r(x)$ counts the number of $x$s in row $x$ and $c(x)$ counts the number of entries in the column of $x$ equal their row number. 

For every element $x\in S$, we have a pair $(r(x),c(x))$ of integers. We can express this data as a polynomial in two variables which we call the {\it shelf polynomial} of $S$:

$$P(S)=\displaystyle\sum_{x\in S}\,t^{r(x)}s^{c(x)}.$$
\end{defn} 

\begin{exmp}

The shelf \vskip 0.1in

\begin{center}
\begin{tabular}{ r| c c c c}
* & 0 & 1 & 2 & 3\\
\hline
  0 & 0 & 1 & 1 & 3\\
  1 & 0 & 1 & 2 & 3\\
  2 & 0 & 1 & 2 & 3\\
  3 & 0 & 1 & 1 & 3\\
\end{tabular} 
\end{center}
\vskip 0.1in
has the shelf polynomial $P(S)=4st$

\end{exmp}

\begin{defn}
A shelf is {\it Latin} or {\it strongly connected} if the shelf operation is left-invertible. This means the rows of a Latin shelf are permutations of $\{1,2,\ldots , n\}$. 
\end{defn} 

We now classify shelves as latin quandles, non-quandle racks, non-rack latin shelves or non-rack shelves. The numbers in parentheses are the numbers assigned to shelves listed in Appendix A. We also present the shelf polynomial of shelves in Appendix A. 

\vskip 0.1in

\textbf{Order 2}

\begin{itemize}
\item Non-rack shelves: (1)
\item Non-quandle Racks: (2)
\end{itemize}

\vskip 0.1in

There are neither non-rack Latin shelves nor Latin quandles of order 2. The shelf polynomials of shelves $(1)$ and $(2)$ are $P(S)=2st$ and $P(S)=2$, respectively. 

\vskip 0.1in

\textbf{Order 3}

\begin{itemize}
\item Non-rack shelves: (5)
\item Non-rack Latin shelves: (1), (2)
\item Non-quandle Racks: (4)
\item Latin Quandles: (3)
\end{itemize}

\vskip 0.1in

The shelf polynomial of the shelf (non-quandle rack) (4) is the constant polynomial $P(S)=3$, and the shelf polynomial of all other shelves is $P(S)=3st$. 

\vskip 0.1in

\textbf{Order 4}

\begin{itemize}
\item Non-rack shelves: (6), (8) -- (12), (14) -- (18)
\item Non-rack Latin shelves: (1), (3), (4), (7), (13)
\item Non-quandle Racks: (5)
\item Latin Quandles: (2) 
\end{itemize}

\vskip 0.1in

The shelf polynomial of the shelves (5) and (6) is the constant polynomial $P(S)=4$, and the shelf polynomial of all other shelves is $P(S)=4st$. 

\vskip 0.1in 

We note to the reader that if $P(S)=|S|st$, then $S$ is clearly idempotent. Then the natural question to ask is ``Is the converse true?''. This leads to the following conjecture strongly supported by our computational results. 

\begin{conj}\label{C2}
If a connected shelf $S$ is idempotent, then $P(S)=|S|st$, where $|S|$ is the cardinality of $S$. 
\end{conj} 

We note to the reader that an idempotent shelf does not always have the shelf polynomial $P(S)=|S|st$. 

For example, the order 3 idempotent shelf, which is non-connected,

\begin{center}
\begin{tabular}{ r| c c c}
* & 0 & 1 & 2\\
\hline
  0 & 0 & 1 & 1\\
  1 & 0 & 1 & 0 \\
  2 & 0 & 1 & 2 \\
\end{tabular} 
\end{center}

has the shelf polynomial $P(S)=3st$, whereas the order 3 idempotent shelf, which is non-connected,

\begin{center}
\begin{tabular}{ r| c c c}
* & 0 & 1 & 2\\
\hline
  0 & 0 & 0 & 0\\
  1 & 1 & 1 & 1 \\
  2 & 2 & 2 & 2 \\
\end{tabular} 
\end{center}

has the shelf polynomial $P(S)=3s^3t^3$. The conjecture says that the shelf polynomial of a connected spindle $S$ is $P(S)=|S|st$, where $|S|$ is the cardinality of $S$.

\section*{Appendix A} \label{Appendix}

In Appendix A, we list all connected shelves of order less than or equal to 5 up to isomorphism. For a full list of shelves of order less than or equal to 5 up to isomorphism and a full list of connected shelves of order less than or equal to six up to isomorphism, we refer the reader to  \cite{Aloglink}. 

\vskip 0.1in

\textbf{Order 2:} There are two connected shelves. 

\begin{tiny}
\begin{enumerate} 
\begin{multicols}{2}
\item $[[0,1],[0,1]]$\vskip 0.2in 
\item $[[1,1],[0,0]]$
\end{multicols}
\end{enumerate} 
\end{tiny} 

\textbf{Order 3:} There are five connected shelves. 

\begin{tiny}
\begin{enumerate} 
\begin{multicols}{2}

\item $[[0,1,2],[0,1,2],[0,1,2]]$ \vskip 0.1in 

\item $[[0,1,2],[0,1,2],[1,0,2]]$ \vskip 0.1in 

\item $[[0,2,1],[2,1,0],[1,0,2]]$ \vskip 0.1in 

\item $[[1,1,1],[2,2,2],[0,0,0]]$ \vskip 0.1in 

\item $[[0,1,1],[0,1,2],[0,1,2]]$ 

\end{multicols}
\end{enumerate} 
\end{tiny} 

\textbf{Order 4:} There are 18 connected shelves.

\begin{tiny}
\begin{enumerate} 
\begin{multicols}{2}

\item $[[0,1,2,3],[0,1,2,3],[0,1,2,3],[0,1,2,3]]$ \vskip 0.1in

\item $[[0,2,3,1],[3,1,0,2],[1,3,2,0],[2,0,1,3]]$ \vskip 0.1in

\item $[[0,1,3,2],[0,1,3,2],[1,0,2,3],[1,0,2,3]]$ \vskip 0.1in

\item $[[0,1,2,3],[0,1,2,3],[1,0,2,3],[1,0,2,3]]$ \vskip 0.1in

\item $[[1,1,1,1],[2,2,2,2],[3,3,3,3],[0,0,0,0]]$ \vskip 0.1in

\item $[[1,1,2,2],[0,0,3,3],[0,0,3,3],[1,1,2,2]]$ \vskip 0.1in

\item $[[0,1,2,3],[0,1,2,3],[0,1,2,3],[1,2,0,3]]$ \vskip 0.1in

\item $[[0,1,1,1],[0,1,2,3],[0,1,2,3],[0,1,2,3]]$ \vskip 0.1in

\item $[[0,1,1,3],[0,1,2,0],[0,1,2,0],[0,1,1,3]]$ \vskip 0.1in

\item $[[0,1,1,3],[0,1,2,3],[0,1,2,3],[0,1,1,3]]$ \vskip 0.1in

\item $[[0,1,1,3],[0,1,2,3],[0,1,2,3],[0,2,2,3]]$ \vskip 0.1in

\item $[[0,1,1,3],[3,1,2,0],[3,1,2,0],[0,1,1,3]]$ \vskip 0.1in

\item $[[0,1,2,3],[0,1,2,3],[0,1,2,3],[0,2,1,3]]$ \vskip 0.1in

\item $[[0,1,1,1],[0,1,2,2],[0,1,2,3],[0,1,2,3]]$ \vskip 0.1in

\item $[[0,1,1,2],[0,1,2,3],[0,1,2,3],[0,1,2,3]]$ \vskip 0.1in

\item $[[0,1,1,3],[0,1,2,3],[0,1,2,3],[0,1,2,3]]$ \vskip 0.1in

\item $[[0,1,1,3],[0,1,2,3],[0,1,2,3],[0,2,1,3]]$ \vskip 0.1in

\item $[[0,1,2,3],[0,1,2,3],[0,1,2,3],[1,0,0,3]]$

\end{multicols}
\end{enumerate} 
\end{tiny}

\textbf{Order 5:} There are 165 connected shelves.

\begin{tiny}
\begin{enumerate} 

\item $[[0,1,1,1,1],[0,1,2,2,2],[0,1,2,3,3],[0,1,2,3,4],[0,1,2,3,4]]$ \vskip 0.1in

\item $[[0,1,1,1,1],[0,1,2,2,2],[0,1,2,3,4],[0,1,2,3,4],[0,1,2,3,4]]$ \vskip 0.1in

\item $[[0,1,1,1,1],[0,1,2,2,4],[0,1,2,3,4],[0,1,2,3,4],[0,1,2,2,4]]$ \vskip 0.1in

\item $[[0,1,1,1,1],[0,1,2,2,4],[0,1,2,3,4],[0,1,2,3,4],[0,1,2,3,4]]$ \vskip 0.1in

\item $[[0,1,1,1,1],[0,1,2,2,4],[0,1,2,3,4],[0,1,2,3,4],[0,1,3,2,4]]$ \vskip 0.1in

\item $[[0,1,1,1,1],[0,1,2,2,4],[0,1,2,3,4],[0,1,2,3,4],[0,1,3,3,4]]$ \vskip 0.1in

\item $[[0,1,1,1,1],[0,1,2,3,4],[0,1,2,3,4],[0,1,2,3,4],[0,1,2,3,4]]$ \vskip 0.1in

\item $[[0,1,1,1,2],[0,1,2,3,3],[0,1,2,3,3],[0,1,2,3,4],[0,1,2,3,4]]$ \vskip 0.1in

\item $[[0,1,1,1,2],[0,1,2,3,4],[0,1,2,3,4],[0,1,2,3,4],[0,1,2,3,4]]$ \vskip 0.1in

\item $[[0,1,1,1,2],[0,1,2,4,4],[0,1,2,4,4],[0,1,2,3,4],[0,1,2,3,4]]$ \vskip 0.1in

\item $[[0,1,1,1,4],[0,1,2,2,0],[0,1,2,3,0],[0,1,2,3,0],[0,1,1,1,4]]$ \vskip 0.1in

\item $[[0,1,1,1,4],[0,1,2,2,4],[0,1,2,3,4],[0,1,2,3,4],[0,1,1,1,4]]$ \vskip 0.1in

\item $[[0,1,1,1,4],[0,1,2,2,4],[0,1,2,3,4],[0,1,2,3,4],[0,1,2,2,4]]$ \vskip 0.1in

\item $[[0,1,1,1,4],[0,1,2,3,0],[0,1,2,3,0],[0,1,2,3,0],[0,1,1,1,4]]$ \vskip 0.1in

\item $[[0,1,1,1,4],[0,1,2,3,4],[0,1,2,3,4],[0,1,2,3,4],[0,1,1,1,4]]$ \vskip 0.1in

\item $[[0,1,1,1,4],[0,1,2,3,4],[0,1,2,3,4],[0,1,2,3,4],[0,1,1,2,4]]$ \vskip 0.1in

\item $[[0,1,1,1,4],[0,1,2,3,4],[0,1,2,3,4],[0,1,2,3,4],[0,1,1,3,4]]$ \vskip 0.1in

\item $[[0,1,1,1,4],[0,1,2,3,4],[0,1,2,3,4],[0,1,2,3,4],[0,1,2,2,4]]$ \vskip 0.1in

\item $[[0,1,1,1,4],[0,1,2,3,4],[0,1,2,3,4],[0,1,2,3,4],[0,1,2,3,4]]$  
\vskip 0.1in

\item $[[0,1,1,1,4],[0,1,2,3,4],[0,1,2,3,4],[0,1,2,3,4],[0,1,3,2,4]]$ \vskip 0.1in

\item $[[0,1,1,1,4],[0,1,2,3,4],[0,1,2,3,4],[0,1,2,3,4],[0,2,1,1,4]]$ \vskip 0.1in

\item $[[0,1,1,1,4],[0,1,2,3,4],[0,1,2,3,4],[0,1,2,3,4],[0,2,1,2,4]]$ \vskip 0.1in

\item $[[0,1,1,1,4],[0,1,2,3,4],[0,1,2,3,4],[0,1,2,3,4],[0,2,1,3,4]]$ \vskip 0.1in

\item $[[0,1,1,1,4],[0,1,2,3,4],[0,1,2,3,4],[0,1,2,3,4],[0,2,2,1,4]]$ 
\vskip 0.1in

\item $[[0,1,1,1,4],[0,1,2,3,4],[0,1,2,3,4],[0,1,2,3,4],[0,2,2,2,4]]$ \vskip 0.1in

\item $[[0,1,1,1,4],[0,1,2,3,4],[0,1,2,3,4],[0,1,2,3,4],[0,2,2,3,4]]$ \vskip 0.1in

\item $[[0,1,1,1,4],[0,1,2,3,4],[0,1,2,3,4],[0,1,2,3,4],[0,2,3,1,4]]$ \vskip 0.1in

\item $[[0,1,1,1,4],[0,1,2,3,4],[0,1,2,3,4],[0,1,2,3,4],[0,2,3,2,4]]$ \vskip 0.1in

\item $[[0,1,1,1,4],[0,1,2,3,4],[0,1,2,3,4],[0,1,2,3,4],[0,2,3,3,4]]$ \vskip 0.1in

\item $[[0,1,1,1,4],[4,1,2,2,0],[4,1,2,3,0],[4,1,2,3,0],[0,1,1,1,4]]$ 

\begin{landscape}
\begin{multicols}{2}

\item $[[0,1,1,1,4],[4,1,2,3,0],[4,1,2,3,0],[4,1,2,3,0],[0,1,1,1,4]]$ \vskip 0.1in

\item $[[0,1,1,2,2],[0,1,2,3,3],[0,1,2,3,3],[0,1,2,3,4],[0,1,2,3,4]]$ \vskip 0.1in

\item $[[0,1,1,2,2],[0,1,2,3,4],[0,1,2,3,4],[0,1,2,3,4],[0,1,2,3,4]]$ \vskip 0.1in

\item $[[0,1,1,2,3],[0,1,2,3,4],[0,1,2,3,4],[0,1,2,3,4],[0,1,2,3,4]]$ \vskip 0.1in

\item $[[0,1,1,2,4],[0,1,2,3,0],[0,1,2,3,0],[0,1,2,3,0],[0,1,1,2,4]]$ \vskip 0.1in

\item $[[0,1,1,2,4],[0,1,2,3,4],[0,1,2,3,4],[0,1,2,3,4],[0,1,1,2,4]]$ \vskip 0.1in

\item $[[0,1,1,2,4],[0,1,2,3,4],[0,1,2,3,4],[0,1,2,3,4],[0,1,1,3,4]]$ \vskip 0.1in

\item $[[0,1,1,2,4],[0,1,2,3,4],[0,1,2,3,4],[0,1,2,3,4],[0,1,2,1,4]]$ \vskip 0.1in

\item $[[0,1,1,2,4],[0,1,2,3,4],[0,1,2,3,4],[0,1,2,3,4],[0,1,2,2,4]]$ \vskip 0.1in

\item $[[0,1,1,2,4],[0,1,2,3,4],[0,1,2,3,4],[0,1,2,3,4],[0,1,2,3,4]]$ \vskip 0.1in

\item $[[0,1,1,2,4],[0,1,2,3,4],[0,1,2,3,4],[0,1,2,3,4],[0,1,3,1,4]]$ \vskip 0.1in

\item $[[0,1,1,2,4],[0,1,2,3,4],[0,1,2,3,4],[0,1,2,3,4],[0,1,3,2,4]]$ \vskip 0.1in

\item $[[0,1,1,2,4],[0,1,2,3,4],[0,1,2,3,4],[0,1,2,3,4],[0,1,3,3,4]]$ \vskip 0.1in

\item $[[0,1,1,2,4],[0,1,2,3,4],[0,1,2,3,4],[0,1,2,3,4],[0,2,1,1,4]]$ \vskip 0.1in

\item $[[0,1,1,2,4],[0,1,2,3,4],[0,1,2,3,4],[0,1,2,3,4],[0,2,1,2,4]]$ 
\vskip 0.1in

\item $[[0,1,1,2,4],[0,1,2,3,4],[0,1,2,3,4],[0,1,2,3,4],[0,2,1,3,4]]$ \vskip 0.1in

\item $[[0,1,1,2,4],[0,1,2,3,4],[0,1,2,3,4],[0,1,2,3,4],[0,2,2,1,4]]$ \vskip 0.1in

\item $[[0,1,1,2,4],[0,1,2,3,4],[0,1,2,3,4],[0,1,2,3,4],[0,2,2,3,4]]$ \vskip 0.1in

\item $[[0,1,1,2,4],[0,1,2,3,4],[0,1,2,3,4],[0,1,2,3,4],[0,2,3,1,4]]$ \vskip 0.1in

\item $[[0,1,1,2,4],[0,1,2,3,4],[0,1,2,3,4],[0,1,2,3,4],[0,2,3,2,4]]$ \vskip 0.1in

\item $[[0,1,1,2,4],[0,1,2,3,4],[0,1,2,3,4],[0,1,2,3,4],[0,2,3,3,4]]$ \vskip 0.1in

\item $[[0,1,1,2,4],[0,1,2,3,4],[0,1,2,3,4],[0,1,2,3,4],[0,3,1,1,4]]$ \vskip 0.1in

\item $[[0,1,1,2,4],[0,1,2,3,4],[0,1,2,3,4],[0,1,2,3,4],[0,3,1,2,4]]$ \vskip 0.1in

\item $[[0,1,1,2,4],[0,1,2,3,4],[0,1,2,3,4],[0,1,2,3,4],[0,3,1,3,4]]$ \vskip 0.1in

\item $[[0,1,1,2,4],[0,1,2,3,4],[0,1,2,3,4],[0,1,2,3,4],[0,3,2,1,4]]$ \vskip 0.1in

\item $[[0,1,1,2,4],[0,1,2,3,4],[0,1,2,3,4],[0,1,2,3,4],[0,3,2,3,4]]$ \vskip 0.1in

\item $[[0,1,1,2,4],[0,1,2,3,4],[0,1,2,3,4],[0,1,2,3,4],[0,3,3,1,4]]$ \vskip 0.1in

\item $[[0,1,1,2,4],[0,1,2,3,4],[0,1,2,3,4],[0,1,2,3,4],[0,3,3,2,4]]$ \vskip 0.1in

\item $[[0,1,1,2,4],[4,1,2,3,0],[4,1,2,3,0],[4,1,2,3,0],[0,1,1,2,4]]$ \vskip 0.1in

\item $[[0,1,1,3,3],[0,1,2,3,2],[0,1,2,3,4],[0,1,2,3,2],[0,1,2,3,4]]$ \vskip 0.1in

\item $[[0,1,1,3,3],[0,1,2,3,3],[0,1,2,3,3],[0,1,1,3,4],[0,1,1,3,4]]$ \vskip 0.1in

\item $[[0,1,1,3,3],[0,1,2,3,3],[0,1,2,3,3],[0,1,2,3,4],[0,1,2,3,4]]$ \vskip 0.1in

\item $[[0,1,1,3,3],[0,1,2,3,3],[0,1,2,3,3],[0,2,1,3,4],[0,2,1,3,4]]$ \vskip 0.1in

\item $[[0,1,1,3,3],[0,1,2,3,3],[0,1,2,3,3],[0,2,2,3,4],[0,2,2,3,4]]$ 
\vskip 0.1in

\item $[[0,1,1,3,3],[0,1,2,3,4],[0,1,2,3,4],[0,1,2,3,4],[0,1,2,3,4]]$ \vskip 0.1in

\item $[[0,1,1,3,3],[0,1,2,3,4],[0,1,2,3,4],[0,2,1,3,4],[0,2,1,3,4]]$ \vskip 0.1in

\item $[[0,1,1,3,3],[0,1,2,3,4],[0,1,2,3,4],[0,2,2,3,4],[0,2,2,3,4]]$ 
\vskip 0.1in

\item $[[0,1,1,3,3],[0,1,2,4,3],[0,1,2,4,3],[0,2,1,3,4],[0,2,1,3,4]]$ \vskip 0.1in

\item $[[0,1,1,3,3],[0,1,2,4,3],[0,1,2,4,3],[0,2,2,3,4],[0,2,2,3,4]]$ \vskip 0.1in

\item $[[0,1,1,3,3],[0,1,2,4,4],[0,1,2,4,4],[0,2,2,3,4],[0,2,2,3,4]]$ \vskip 0.1in

\item $[[0,1,1,3,4],[0,1,2,0,4],[0,1,2,0,4],[0,1,1,3,4],[0,1,1,3,4]]$ \vskip 0.1in

\item $[[0,1,1,3,4],[0,1,2,0,4],[0,1,2,0,4],[0,1,1,3,4],[0,1,2,3,4]]$ \vskip 0.1in

\item $[[0,1,1,3,4],[0,1,2,0,4],[0,1,2,0,4],[0,1,1,3,4],[0,2,1,0,4]]$ \vskip 0.1in

\item $[[0,1,1,3,4],[0,1,2,0,4],[0,1,2,0,4],[0,1,1,3,4],[0,2,1,3,4]]$ \vskip 0.1in

\item $[[0,1,1,3,4],[0,1,2,0,4],[0,1,2,0,4],[0,1,1,3,4],[0,2,2,3,4]]$ \vskip 0.1in

\item $[[0,1,1,3,4],[0,1,2,0,4],[0,1,2,0,4],[0,1,1,3,4],[3,2,1,0,4]]$ \vskip 0.1in

\item $[[0,1,1,3,4],[0,1,2,0,4],[0,1,2,0,4],[0,1,1,3,4],[3,2,1,3,4]]$ \vskip 0.1in

\item $[[0,1,1,3,4],[0,1,2,3,3],[0,1,2,3,3],[0,1,2,3,4],[0,1,2,3,4]]$ \vskip 0.1in

\item $[[0,1,1,3,4],[0,1,2,3,3],[0,1,2,3,3],[0,2,1,3,4],[0,2,1,3,4]]$ \vskip 0.1in

\item $[[0,1,1,3,4],[0,1,2,3,4],[0,1,2,3,4],[0,1,1,3,4],[0,1,1,3,4]]$ \vskip 0.1in

\item $[[0,1,1,3,4],[0,1,2,3,4],[0,1,2,3,4],[0,1,1,3,4],[0,1,2,3,4]]$ \vskip 0.1in

\item $[[0,1,1,3,4],[0,1,2,3,4],[0,1,2,3,4],[0,1,1,3,4],[0,2,1,0,4]]$ \vskip 0.1in

\item $[[0,1,1,3,4],[0,1,2,3,4],[0,1,2,3,4],[0,1,1,3,4],[0,2,1,3,4]]$ \vskip 0.1in

\item $[[0,1,1,3,4],[0,1,2,3,4],[0,1,2,3,4],[0,1,1,3,4],[0,2,2,3,4]]$ \vskip 0.1in

\item $[[0,1,1,3,4],[0,1,2,3,4],[0,1,2,3,4],[0,1,1,3,4],[3,0,0,0,4]]$ \vskip 0.1in

\item $[[0,1,1,3,4],[0,1,2,3,4],[0,1,2,3,4],[0,1,1,3,4],[3,0,3,0,4]]$ \vskip 0.1in

\item $[[0,1,1,3,4],[0,1,2,3,4],[0,1,2,3,4],[0,1,1,3,4],[3,1,1,0,4]]$ \vskip 0.1in

\item $[[0,1,1,3,4],[0,1,2,3,4],[0,1,2,3,4],[0,1,1,3,4],[3,1,2,0,4]]$ \vskip 0.1in

\item $[[0,1,1,3,4],[0,1,2,3,4],[0,1,2,3,4],[0,1,1,3,4],[3,2,1,0,4]]$ \vskip 0.1in

\item $[[0,1,1,3,4],[0,1,2,3,4],[0,1,2,3,4],[0,1,1,3,4],[3,2,2,0,4]]$ \vskip 0.1in

\item $[[0,1,1,3,4],[0,1,2,3,4],[0,1,2,3,4],[0,1,2,3,1],[0,1,2,3,4]]$ \vskip 0.1in

\item $[[0,1,1,3,4],[0,1,2,3,4],[0,1,2,3,4],[0,1,2,3,2],[0,1,2,3,4]]$ \vskip 0.1in

\item $[[0,1,1,3,4],[0,1,2,3,4],[0,1,2,3,4],[0,1,2,3,4],[0,1,2,3,4]]$ \vskip 0.1in

\item $[[0,1,1,3,4],[0,1,2,3,4],[0,1,2,3,4],[0,1,2,3,4],[0,1,3,2,4]]$ \vskip 0.1in

\item $[[0,1,1,3,4],[0,1,2,3,4],[0,1,2,3,4],[0,1,2,3,4],[0,1,3,3,4]]$ \vskip 0.1in

\item $[[0,1,1,3,4],[0,1,2,3,4],[0,1,2,3,4],[0,1,2,3,4],[0,2,1,1,4]]$ \vskip 0.1in

\item $[[0,1,1,3,4],[0,1,2,3,4],[0,1,2,3,4],[0,1,2,3,4],[0,2,1,2,4]]$ 
\vskip 0.1in

\item $[[0,1,1,3,4],[0,1,2,3,4],[0,1,2,3,4],[0,1,2,3,4],[0,2,1,3,4]]$ 
\vskip 0.1in

\item $[[0,1,1,3,4],[0,1,2,3,4],[0,1,2,3,4],[0,1,2,3,4],[0,2,2,3,4]]$ \vskip 0.1in

\item $[[0,1,1,3,4],[0,1,2,3,4],[0,1,2,3,4],[0,1,2,3,4],[0,2,3,1,4]]$ \vskip 0.1in

\item $[[0,1,1,3,4],[0,1,2,3,4],[0,1,2,3,4],[0,1,2,3,4],[0,2,3,2,4]]$ \vskip 0.1in

\item $[[0,1,1,3,4],[0,1,2,3,4],[0,1,2,3,4],[0,1,2,3,4],[0,3,1,1,4]]$ \vskip 0.1in

\item $[[0,1,1,3,4],[0,1,2,3,4],[0,1,2,3,4],[0,1,2,3,4],[0,3,1,2,4]]$ \vskip 0.1in

\item $[[0,1,1,3,4],[0,1,2,3,4],[0,1,2,3,4],[0,1,2,3,4],[0,3,2,1,4]]$ 
\vskip 0.1in

\item $[[0,1,1,3,4],[0,1,2,3,4],[0,1,2,3,4],[0,1,2,3,4],[0,3,3,1,4]]$ \vskip 0.1in

\item $[[0,1,1,3,4],[0,1,2,3,4],[0,1,2,3,4],[0,1,2,3,4],[0,3,3,2,4]]$ \vskip 0.1in

\item $[[0,1,1,3,4],[0,1,2,3,4],[0,1,2,3,4],[0,2,1,3,4],[0,2,1,3,4]]$ \vskip 0.1in

\item $[[0,1,1,3,4],[0,1,2,3,4],[0,1,2,3,4],[0,2,1,3,4],[0,2,2,3,4]]$ \vskip 0.1in

\item $[[0,1,1,3,4],[0,1,2,4,3],[0,1,2,4,3],[0,1,1,3,4],[0,1,1,3,4]]$ \vskip 0.1in

\item $[[0,1,1,3,4],[0,1,2,4,3],[0,1,2,4,3],[0,1,2,3,4],[0,1,2,3,4]]$ \vskip 0.1in

\item $[[0,1,1,3,4],[0,1,2,4,3],[0,1,2,4,3],[0,2,1,3,4],[0,2,1,3,4]]$ \vskip 0.1in

\item $[[0,1,1,3,4],[0,1,2,4,3],[0,1,2,4,3],[0,2,2,3,4],[0,2,2,3,4]]$ \vskip 0.1in

\item $[[0,1,1,3,4],[3,1,2,0,0],[3,1,2,0,0],[0,1,1,3,4],[0,1,1,3,4]]$ \vskip 0.1in

\item $[[0,1,1,3,4],[3,1,2,0,4],[3,1,2,0,4],[0,1,1,3,4],[0,1,2,3,4]]$ \vskip 0.1in

\item $[[0,1,1,3,4],[3,1,2,0,4],[3,1,2,0,4],[0,1,1,3,4],[0,2,1,0,4]]$ \vskip 0.1in

\item $[[0,1,1,3,4],[3,1,2,0,4],[3,1,2,0,4],[0,1,1,3,4],[0,2,1,3,4]]$ \vskip 0.1in

\item $[[0,1,1,3,4],[3,1,2,0,4],[3,1,2,0,4],[0,1,1,3,4],[3,1,1,0,4]]$ \vskip 0.1in

\item $[[0,1,1,3,4],[3,1,2,0,4],[3,1,2,0,4],[0,1,1,3,4],[3,1,2,0,4]]$ \vskip 0.1in

\item $[[0,1,1,3,4],[3,1,2,0,4],[3,1,2,0,4],[0,1,1,3,4],[3,2,1,0,4]]$ \vskip 0.1in

\item $[[0,1,1,3,4],[3,1,2,0,4],[3,1,2,0,4],[0,1,1,3,4],[3,2,2,0,4]]$ \vskip 0.1in

\item $[[0,1,1,3,4],[3,1,2,4,0],[3,1,2,4,0],[0,1,1,3,4],[0,1,1,3,4]]$ \vskip 0.1in

\item $[[0,1,1,4,3],[0,1,2,3,4],[0,1,2,3,4],[0,1,2,3,4],[0,1,2,3,4]]$ \vskip 0.1in

\item $[[0,1,1,4,3],[0,1,2,3,4],[0,1,2,3,4],[0,2,1,3,4],[0,2,1,3,4]]$ 
\vskip 0.1in

\item $[[0,1,1,4,3],[0,1,2,4,3],[0,1,2,4,3],[0,1,2,3,4],[0,1,2,3,4]]$
\vskip 0.1in

\item $[[0,1,1,4,3],[0,1,2,4,3],[0,1,2,4,3],[0,2,1,3,4],[0,2,1,3,4]]$ 
\vskip 0.1in

\item $[[0,1,2,3,4],[0,1,2,3,4],[0,1,2,3,4],[0,1,2,3,4],[0,1,2,3,4]]$ \vskip 0.1in

\item $[[0,1,2,3,4],[0,1,2,3,4],[0,1,2,3,4],[0,1,2,3,4],[0,1,3,2,4]]$ \vskip 0.1in

\item $[[0,1,2,3,4],[0,1,2,3,4],[0,1,2,3,4],[0,1,2,3,4],[0,2,1,1,4]]$ \vskip 0.1in

\item $[[0,1,2,3,4],[0,1,2,3,4],[0,1,2,3,4],[0,1,2,3,4],[0,2,3,1,4]]$ \vskip 0.1in

\item $[[0,1,2,3,4],[0,1,2,3,4],[0,1,2,3,4],[0,1,2,3,4],[1,0,0,0,4]]$ \vskip 0.1in

\item $[[0,1,2,3,4],[0,1,2,3,4],[0,1,2,3,4],[0,1,2,3,4],[1,0,0,1,4]]$ \vskip 0.1in

\item $[[0,1,2,3,4],[0,1,2,3,4],[0,1,2,3,4],[0,1,2,3,4],[1,0,0,2,4]]$ \vskip 0.1in

\item $[[0,1,2,3,4],[0,1,2,3,4],[0,1,2,3,4],[0,1,2,3,4],[1,0,3,2,4]]$ \vskip 0.1in

\item $[[0,1,2,3,4],[0,1,2,3,4],[0,1,2,3,4],[0,1,2,3,4],[1,2,0,0,4]]$ \vskip 0.1in

\item $[[0,1,2,3,4],[0,1,2,3,4],[0,1,2,3,4],[0,1,2,3,4],[1,2,3,0,4]]$ \vskip 0.1in

\item $[[0,1,2,3,4],[0,1,2,3,4],[0,1,2,3,4],[0,2,1,3,4],[0,2,1,3,4]]$ \vskip 0.1in

\item $[[0,1,2,3,4],[0,1,2,3,4],[0,1,2,3,4],[0,2,1,3,4],[1,0,0,3,4]]$ \vskip 0.1in

\item $[[0,1,2,3,4],[0,1,2,3,4],[0,1,2,3,4],[0,2,1,3,4],[1,0,1,3,4]]$ \vskip 0.1in

\item $[[0,1,2,3,4],[0,1,2,3,4],[0,1,2,3,4],[0,2,1,3,4],[1,0,2,3,4]]$ \vskip 0.1in

\item $[[0,1,2,3,4],[0,1,2,3,4],[0,1,2,3,4],[0,2,1,3,4],[1,2,0,3,4]]$ \vskip 0.1in

\item $[[0,1,2,3,4],[0,1,2,3,4],[0,1,2,3,4],[0,2,1,3,4],[1,2,1,3,4]]$ 
\vskip 0.1in

\item $[[0,1,2,3,4],[0,1,2,3,4],[0,1,2,3,4],[0,2,1,3,4],[3,3,3,0,4]]$ \vskip 0.1in

\item $[[0,1,2,3,4],[0,1,2,3,4],[0,1,2,3,4],[1,0,0,3,4],[1,0,0,3,4]]$ \vskip 0.1in

\item $[[0,1,2,3,4],[0,1,2,3,4],[0,1,2,3,4],[1,0,0,3,4],[1,0,1,3,4]]$ \vskip 0.1in

\item $[[0,1,2,3,4],[0,1,2,3,4],[0,1,2,3,4],[1,0,0,3,4],[1,2,0,3,4]]$ \vskip 0.1in

\item $[[0,1,2,3,4],[0,1,2,3,4],[0,1,2,3,4],[1,0,0,3,4],[1,2,1,3,4]]$ \vskip 0.1in

\item $[[0,1,2,3,4],[0,1,2,3,4],[0,1,2,3,4],[1,0,0,3,4],[2,0,0,3,4]]$ \vskip 0.1in

\item $[[0,1,2,3,4],[0,1,2,3,4],[0,1,2,3,4],[1,0,0,3,4],[2,0,1,3,4]]$ 
\vskip 0.1in

\item $[[0,1,2,3,4],[0,1,2,3,4],[0,1,2,3,4],[1,0,0,3,4],[2,2,1,3,4]]$ \vskip 0.1in

\item $[[0,1,2,3,4],[0,1,2,3,4],[0,1,2,3,4],[1,2,0,3,4],[1,2,0,3,4]]$ 
\vskip 0.1in

\item $[[0,1,2,3,4],[0,1,2,3,4],[0,1,2,3,4],[1,2,0,3,4],[2,0,1,3,4]]$ \vskip 0.1in

\item $[[0,1,2,3,4],[0,1,2,3,4],[0,1,2,4,3],[1,0,2,3,4],[1,0,2,3,4]]$ \vskip 0.1in

\item $[[0,1,2,3,4],[0,1,2,3,4],[1,0,2,3,4],[1,0,2,3,4],[1,0,2,3,4]]$ \vskip 0.1in

\item $[[0,1,2,3,4],[0,1,2,3,4],[1,0,2,3,4],[1,0,2,3,4],[1,0,3,2,4]]$ \vskip 0.1in

\item $[[0,1,2,3,4],[0,1,2,4,3],[0,1,2,4,3],[0,2,1,3,4],[0,2,1,3,4]]$ \vskip 0.1in

\item $[[0,1,2,4,3],[0,1,2,4,3],[0,1,2,4,3],[0,2,1,3,4],[0,2,1,3,4]]$ \vskip 0.1in

\item $[[0,1,2,4,3],[0,1,2,4,3],[0,1,2,4,3],[1,0,0,3,4],[1,0,0,3,4]]$ \vskip 0.1in

\item $[[0,1,2,4,3],[0,1,2,4,3],[0,1,2,4,3],[1,2,0,3,4],[1,2,0,3,4]]$ \vskip 0.1in

\item $[[0,1,2,4,3],[0,1,2,4,3],[1,0,2,4,3],[1,0,2,3,4],[1,0,2,3,4]]$ \vskip 0.1in

\item $[[0,1,2,4,3],[0,1,2,4,3],[1,0,2,4,3],[2,2,0,3,4],[2,2,1,3,4]]$ 
\vskip 0.1in

\item $[[0,1,3,2,2],[0,1,4,2,2],[3,4,2,0,1],[2,2,0,3,4],[2,2,1,3,4]]$ \vskip 0.1in

\item $[[0,2,1,4,3],[3,1,4,0,2],[4,3,2,1,0],[2,4,0,3,1],[1,0,3,2,4]]$ \vskip 0.1in

\item $[[0,2,3,4,1],[2,1,4,0,3],[3,4,2,1,0],[4,0,1,3,2],[1,3,0,2,4]]$ \vskip 0.1in

\item $[[0,2,3,4,1],[3,1,4,2,0],[4,0,2,1,3],[1,4,0,3,2],[2,3,1,0,4]]$ \vskip 0.1in

\item $[[1,1,1,1,1],[2,2,2,2,2],[3,3,3,3,3],[4,4,4,4,4],[0,0,0,0,0]]$ \vskip 0.1in
\end{multicols}{2}
\end{landscape}
\end{enumerate} 
\end{tiny}

\section*{Acknowledgement} 
Mohamed Elhamdadi was partially supported by Simons Foundation collaboration grant 712462. The authors would like to thank Sujoy Mukhejree and Manpreet Singh for fruitfull comments which improved the paper.

\end{document}